\documentclass[12pt]{amsart}

\usepackage{amssymb}
\usepackage{amsmath}
\numberwithin{equation}{section}
\usepackage{breakcites}
\usepackage{color}
\usepackage{graphicx}
\usepackage{epstopdf}
\usepackage{array}
\usepackage{hyperref}

\DeclareMathOperator{\sgn}{sgn}
\DeclareMathOperator{\im}{Im}
\usepackage[vcentermath,enableskew]{youngtab}
\usepackage{ytableau}
\usepackage{amssymb}
\usepackage{rotating}
\usepackage{young}
\usepackage{enumerate}
\usepackage{mathtools}
\usepackage{tikz}
\usetikzlibrary{calc,decorations.pathreplacing,calligraphy,matrix}

\newcommand{\be}{\begin{equation}} 
\newcommand{\bea}{\begin{eqnarray}}
\newcommand{\ee}{\end{equation}}
\newcommand{\beas}{\begin{eqnarray*}}
\newcommand{\eea}{\end{eqnarray}}
\newcommand{\eeas}{\end{eqnarray*}}

\newcommand{\cal}{\mathcal}

\def\C{{\mathbb C}}

\def\z3{{\mathbb Z_3}}

\newcommand{\Tlam}{\mathcal{T}_\lambda}

\def\sg{\mathfrak S}
\newtheorem{theorem}{Theorem}[section]
\newtheorem{definition}[theorem]{Definition}

\newtheorem{corollary}[theorem]{Corollary}

\newtheorem{lem}[theorem]{Lemma}

\begin{document}
\title[The construction of a class  of presentations for Specht modules]{The construction of a class  \\ of presentations for Specht modules}
\author[Friedmann]{Tamar Friedmann}
\address{Department of Mathematics, Colby College}
\email{tfriedma@colby.edu}

\begin{abstract} We build on  the methods introduced by Friedmann, Hanlon, Stanley, and Wachs, and further developed by Brauner and Friedmann, to construct additional classes of presentations of Specht modules. We obtain these presentations by defining a linear operator which is a symmetrized sum of dual Garnir relations on the space of column tabloids. Our presentations apply to the vast majority of shapes of Specht modules. 
\end{abstract}

\maketitle

\section{Introduction}

The Specht modules $S^\lambda$, where $\lambda$ is a partition of $n$, 
give a complete set of  irreducible 
representations of the symmetric group $\sg_n$ over a field  of characteristic $0$, say $\mathbb C $.
They can be constructed as  
subspaces of the regular representation $\mathbb C \sg_n$  or  as presentations given in terms of generators and relations, 
known as Garnir relations. This paper deals primarily with the latter type of construction.

Let $\lambda= (\lambda_1 \geq\dots \geq \lambda_l)$ be a partition of $n$.  A {\em Young tableaux} of shape $\lambda$ is a filling of the Young diagram of shape $\lambda$ with distinct entries from the set $[n]:=\{1,2,\dots,n\}$. Let $\mathcal T_\lambda$ be the set of  Young tableaux of shape
$\lambda$. The symmetric group $\sg_n$ acts on $\mathcal T_\lambda$ by replacing
 each entry of a tableau by its image under the permutation in $\sg_n$.

To construct the Specht module as a submodule of the regular representation, one can use  Young symmetrizers.  For $t \in \mathcal T_\lambda$, the {\em Young symmetrizer} is defined by 
\begin{equation} \label{et} e_t:=\sum_{\alpha \in R_{t}}  \alpha  \\  \sum_{\beta \in C_t} \sgn(\beta) \beta ,\end{equation}
where $C_t$ is the column stabilizer of $t$ and $R_t$ is the row stabilizer of $t$.
The {\em Specht module} $S^{\lambda}$ is the submodule of the regular representation  $\C \sg_{n}$ spanned by $\{\tau e_t : \tau \in \sg_{n} \} $.

To construct the Specht module as a presentation, one can use column tabloids and Garnir relations. Let $M^\lambda$ be the vector space (over $\C$) generated by  
$\mathcal T_\lambda$ subject only to column relations, which are of
the form $t+s$, where $s\in \mathcal T_\lambda$ is obtained from $t\in \mathcal T_\lambda$ by switching two
entries in the same column.  
Given $t \in \mathcal T_\lambda$, let $[t]$ 
denote the coset of $t$
in $M^\lambda$. These cosets, which are called {\it column tabloids}, generate $M^{\lambda}$.
A Young tableau is {\em column strict} if the entries of each of its columns increase from top to bottom.  Clearly,
$\{[t]  :  t \mbox{ is a column strict Young tableau of shape } \lambda\}$ is a basis for $M^\lambda$.

In \cite{Fu}, Fulton introduces dual Garnir relations on the column tabloids and shows that $S^\lambda$ is isomorphic to the quotient space of $M^\lambda$ by these relations. There is a dual Garnir relation for each $t\in \Tlam$, each choice of adjacent columns, and each $\ell$ up to the length of the next column. In particular, each Garnir relation is a weighted sum of a column tabloid $[t]$ and  column tabloids obtained from $t$ by exchanging $\ell$ entries of a column with the top $\ell$ entries of the next column. Fulton then obtains a simplification: it is enough to use only the dual Garnir relations that exchange exactly one entry of a column with the top entry of the next column, i.e. we can just fix $\ell=1$. 

An analogous simplification  is obtained in \cite{FHW3}, which improves upon a result in \cite{FHSW}. In the presentation of \cite{FHW3}, $\ell$ is also restricted to a single value, but this time it is the maximum possible value: $\ell$ equals the length of the next column, so as many entries are exchanged as the shape of $\lambda$ allows. This presentation holds for partitions whose conjugate has distinct parts. 

A different simplification is obtained in \cite{BF}, where a symmetrized sum of dual Garnir relations with $\ell=1$ is introduced. The number of relations needed is dramatically reduced: the construction uses a single relation for every pair of adjacent columns and $[t]$ varies in $M^\lambda$, a significantly smaller space than $\mathcal T_\lambda$. The presentation of \cite{BF} holds for all partitions.

In the present paper we  consider intermediate values of $\ell$. For what shapes $\lambda$ would dual Garnir relations that exchange exactly $\ell$ entries between the columns provide a presentation for $S^\lambda$? This is a question posed in \cite{FHW3}. Our results answer that question and generalize both \cite{FHW3}, where $\ell$ is maximal, and \cite{BF}, where $\ell =1$. Our methods are those envisioned in \cite{FHW3}. The question posed there also  inspired \cite{MMS} to address it using a different approach, via representations of the general linear group. See Remark on p. \pageref{MMSrelation} following Theorem \ref{spechtgarnirell} for a discussion relating the results of \cite{MMS} to the results in the current paper. 

Our main result is contained in Theorems \ref{imetaellkeralpha} and \ref{spechtgarnirell}. In Theorem \ref{imetaellkeralpha}, we provide conditions on the shape of a 2-column partition $\mu =(n,m)'$ for any $\ell$, using eigenvalues of an operator, $\eta_\ell$, on $M^\mu$. Given a value of $\ell$, for any $\mu$ for which the conditions are satisfied, we have obtained a presentation of $S^\mu$. In Theorem \ref{spechtgarnirell}, we use the results for 2-column shapes to state conditions on partitions $\lambda$ with any number of columns. In Table \ref{tab:wli}, we provide some computer-generated data that tells us which 2-column shapes satisfy these conditions. 
The data indicates that the conditions are satisfied for the vast majority of shapes $\mu$ and values  of $\ell$. 

The work in \cite{FHSW} was presented in the language of the generalized Jacobi relations that define  the LAnKe or 
Filippov algebra \cite{BL, DI, DT, Fi, Fr, Gu, Ka, Li, Ta}. An observation in   \cite{FHSW}, that the  restricted class of Garnir relations that fix $\ell$ to be the maximum possible value  (for staircase partitions $\lambda$) corresponds to the generalized Jacobi relations, motivated the  work in the papers \cite{FHSW, FHW3, BF, MMS} and the current paper. 
 
This paper is organized as follows. In Section \ref{notation}, we introduce notation and review the relevant prior results. In Section \ref{eigens}, we introduce the relevant relation  and compute its eigenvalues $w_{\ell, i}$, which appear in Theorem \ref{etaell}. A combinatorial identity emerges, given in Corollary \ref{id}. Section \ref{class} contains our main results, along with a statement about the equivalence of two sets of combinatorial conditions (Corollary \ref{equivcomb}).

\section{Notation and prior results}\label{notation}

The symmetric group $\sg_n$ acts on $\mathcal T_\lambda$ by replacing
 each entry of a tableau by its image under the permutation in $\sg_n$.
 This induces a representation of $\sg_n$ on $M^\lambda$.
In \cite[Ch. 7.4]{Fu}, Fulton introduces a map  $$\alpha: {M}^{\lambda} \to S^{\lambda}$$  
given by 
\[ \alpha: [t] \mapsto e_{t}. \]
The map $\alpha$ is $\sg_{n}$-equivariant and surjective. Moreover, $\ker(\alpha)$ is generated by a set of relations which Fulton calls the \emph{dual Garnir relations}.

The dual Garnir relation $g_{c,\ell}(t)$ is
\begin{equation} \label{garnir} g_{c,\ell}(t) = [t] - \pi_{c,\ell}(t) ,\end{equation}
where $\pi_{c,\ell}(t)$ is the sum of column tabloids obtained from all possible ways of exchanging the top $\ell$ elements of the $(c+1)^{st}$ column of $t$ with any subset  of size $\ell$ of the elements of column $c$, preserving the vertical order of each set of $\ell$ elements. 
Note that $t$ can be any tableau, not necessarily with increasing columns.

 If we let 
$G^{\lambda}$ be the subspace of $M^\lambda$ generated by the Garnir relations in
 \begin{equation} \label{garnireq} \{ g_{c,\ell}(t) : c \in [\lambda_1-1], \ell \in [\lambda'_{c+1}], t \in \mathcal T_\lambda \}  ,\end{equation}  
 where $\lambda '$ is the conjugate partition of $\lambda$, then 
 $G^\lambda$ is invariant under the action of
$\sg_n$. 
In \cite[Ch. 7.4]{Fu}, Fulton shows that $G^\lambda=\ker(\alpha)$, thereby obtaining  the following   presentation of $S^\lambda$:
   \be \label{Fultoneq} M^\lambda/G^\lambda \cong_{\sg_n} S^\lambda. \ee
As mentioned in the introduction, $G^\lambda$ contains a dual Garnir relation for each $t\in \Tlam$, each choice of adjacent columns, and each $\ell$ up to the length of the next column.

 On page 102 (after Ex.~15) of \cite{Fu}, a presentation of $S^\lambda$ with a smaller set of relations is given.  In this presentation, the index $\ell$ in $g_{c,\ell} (t)$ of (\ref{garnireq}) is restricted to a single value:  $\ell=\min [\lambda'_{c+1}]= 1$. More precisely, the presentation is
\be \label{1fulton} M^{\lambda} / G^{\lambda,\min} \cong_{\sg_n} S^{\lambda} ,\ee where
 $G^{\lambda,\min}$ is the subspace of $G^\lambda$ generated by the subset of Garnir relations with $\ell=1$: $$\{ g_{c,1}(t) : c \in [\lambda_1-1], t \in \mathcal T_\lambda \}  .$$ 
 This is Fulton's simplification mentioned in the introduction. 
 
 In the analogous simplification of \cite{FHW3},
  the index $\ell$ in $g_{c,\ell} (t)$ of (\ref{garnireq}) is restricted to the maximum value $\ell= \lambda'_{c+1}$:

\begin{theorem} \cite[Theorem 1.1]{FHW3} \label{fhw3} Let $\lambda$ be a partition whose conjugate has distinct parts. Then 
\begin{equation} \label{maineq} M^{\lambda} / G^{\lambda,\max} \cong_{\sg_n} S^{\lambda} ,\end{equation} where
 $G^{\lambda,\max}$ is the subspace of $G^\lambda$ generated by $$\{ g_{c,\lambda'_{c+1}} (t): c \in [\lambda_1-1], t \in \mathcal T_\lambda \}  .$$
Moreover, this set of relations can be further reduced by restricting $t$ to the set of column strict tableaux.
\end{theorem}

The approach in \cite{FHW3}, which improves on an earlier result \cite{FHSW} that applied only to staircase partitions, is to define a certain linear operator on the space of column tabloids and study its eigenspaces. This approach has also been used in \cite{BF}  to obtain a different presentation of $S^\lambda$ with a reduced number of relations, which  works for all 
shapes. Rather than using a subset of the Garnir relations,  \cite{BF} introduce a relation consisting of symmetrized sums of  the dual Garnir relations that 
generate $G^{\lambda, \min}$: 
\begin{equation} \label{BFsym} \eta _{c,1}([t]) = m[t] - \sum [s] \; ,\end{equation}
where the sum ranges over all possible tableaux $s$ obtained from $t$ by swapping one entry in column $c+1$ of $t$ with one entry in column $c$.

The relation $\eta_{c,1}$ can be thought of as a sum of $g_{c,1}$ relations which has the advantage of symmetrizing over all positions of elements in the $(c+1)^{st}$ column. 
For $t \in \mathcal{T}_{\lambda}$, let $h_{c,1}([t])$ be the image of $\eta_{c,1}$ on the $c$ and $(c+1)^{st}$ columns of $[t]$ that leaves the other columns of $[t]$ fixed.  
The result obtained in \cite{BF} is:

\begin{theorem}\cite[Theorem 3.5]{BF}\label{bf}
For any partition $\lambda$ of $n$, let ${H}^{\lambda}$ be the space generated by $h_{c,1}([t])$ for every $[t] \in M^\lambda$ and $1 \leq c \leq \lambda_{1}-1$. 
Then the kernel of $\alpha$ is ${H}^{\lambda}$. Thus,
\[  {M}^{\lambda} / {H}^{\lambda} \cong_{\sg_n} S^{\lambda}  .\]
\end{theorem}
We can see that only the single relation $\eta _{c,1}$ is needed for each pair of adjacent columns. 

In the next section, we address the case of an intermediate value of $\ell$ that was proposed in \cite{FHW3}, generalizing the methods of \cite{FHW3,BF}. As mentioned in the introduction, \cite{MMS} was motivated by \cite{FHW3} to address this same question in a different way, using representations of the general linear group. See Remark on p. \pageref{MMSrelation} following Theorem \ref{spechtgarnirell}  for a statement of their theorem and a discussion relating their results to ours. 

\section{A linear operator and its eigenvalues}\label{eigens}
We now define a linear operator on the space of column tabloids, and study its eigenspaces. 
We use a symmetrized sum of $g_{c,\ell}(t)$ relations, where we symmetrize over all ways of exchanging any $\ell$ elements of the $(c+1)^{st}$ column with any $\ell$ elements of column $c$, where $\ell\in [\lambda '_{c+1}]$. This is a generalization of \cite{BF}, where $\ell=1$, and of \cite{FHW3}, where $\ell =\lambda'_{c+1}$.

In this section we focus on 2-column partitions $\mu$ of $n+m$ with shape $2^{m}1^{n-m}$, so $\mu$ has a column of size $n$ and a column of size $m$ for $1 \leq m \leq n$, and $\mu'=(n,m)$. We shall address the implications of these results to partitions with more than two columns in the next section. 

\begin{definition} \label{defeta}
Let $\mu = 2^m1^{n-m}$,  and let $\ell\in [m]$. We define the map  $\eta_\ell: {M}^{\mu} \to {M}^{\mu}$ to be 
\[\eta_\ell [t] = {m\choose \ell}[t] - \sum [s] \]
where the sum ranges over all possible tableaux $s$ obtained from $t$ by swapping $\ell$ entries in the second column of $t$ with $\ell$ entries in the first column, preserving the vertical order of each set of $\ell$ entries. 
\end{definition}

\noindent {\bf Example:} Let
\small
\[ [t] =  \hspace{.25em} \begin{array}{|c|c|}
  1 & 5\\
  2 & 6 \\
  3 &7\\
  4\\
\end{array}. \]
\normalsize
Then
\small
\begin{align*} \eta_2([t]) =  3 \hspace{.5em}  \begin{array}{|c|c|}
  1 & 5\\
  2 & 6 \\
  3 &7\\
  4\\
  \end{array}
   \hspace{.25em} &- \hspace{.25em}
  \left( \hspace{.25em} 
  \begin{array}{|c|c|}
  \color{red} 5 & \color{red} 1\\
  \color{red} 6 & \color{red}2 \\
  3&7\\
  4\\
\end{array}\hspace{.25em} + \hspace{.25em} 
 \begin{array}{|c|c|}
  \color{red} 5 & \color{red} 1\\
  2 & \color{red}3 \\
  \color{red}6&7\\
  4\\
\end{array}\hspace{.25em} + \hspace{.25em} 
 \begin{array}{|c|c|}
  \color{red} 5 & \color{red}1\\
  2& \color{red}4 \\
  3&7\\
  \color{red} 6 \\
\end{array}\hspace{.25em} + \hspace{.25em} 
\begin{array}{|c|c|}
  1 & \color{red}2\\
  \color{red} 5& \color{red}3 \\
  \color{red} 6&7\\
  4 \\
\end{array}\hspace{.25em} + \hspace{.25em} 
\begin{array}{|c|c|}
  1 & \color{red}2\\
  \color{red} 5& \color{red}4 \\
  3&7\\
  \color{red} 6 \\
\end{array}\hspace{.25em} + \hspace{.25em} 
\begin{array}{|c|c|}
  1 & \color{red}3\\
   2& \color{red}4 \\
  \color{red}5&7\\
  \color{red} 6 \\
\end{array}\hspace{.25em} 
 \hspace{.25em}   \right) 
 \\ 
 &-\hspace{.25em}
  \left( \hspace{.25em} 
  \begin{array}{|c|c|}
  \color{red} 5 & \color{red} 1\\
  \color{red} 7 & 6 \\
  3&\color{red}2\\
  4\\
\end{array}\hspace{.25em} + \hspace{.25em} 
 \begin{array}{|c|c|}
  \color{red} 5 & \color{red} 1\\
  2 & 6 \\
  \color{red}7&\color{red}3\\
  4\\
\end{array}\hspace{.25em} + \hspace{.25em} 
 \begin{array}{|c|c|}
  \color{red} 5 & \color{red}1\\
  2& 6 \\
  3&\color{red}4\\
  \color{red} 7 \\
\end{array}\hspace{.25em} + \hspace{.25em} 
\begin{array}{|c|c|}
  1 & \color{red}2\\
  \color{red} 5& 6 \\
  \color{red} 7&\color{red}3\\
  4 \\
\end{array}\hspace{.25em} + \hspace{.25em} 
\begin{array}{|c|c|}
  1 & \color{red}2\\
  \color{red} 5& 6 \\
  3&\color{red}4\\
  \color{red} 7 \\
\end{array}\hspace{.25em} + \hspace{.25em} 
\begin{array}{|c|c|}
  1 & \color{red}3\\
   2& 6 \\
  \color{red}5&\color{red}4\\
  \color{red} 7 \\
\end{array}\hspace{.25em} 
 \hspace{.25em}   \right)
 \\
 &-\hspace{.25em}
  \left( \hspace{.25em} 
  \begin{array}{|c|c|}
  \color{red} 6 & 5\\
  \color{red} 7 & \color{red} 1 \\
  3&\color{red}2\\
  4\\
\end{array}\hspace{.25em} + \hspace{.25em} 
 \begin{array}{|c|c|}
  \color{red} 6 & 5\\
  2 & \color{red} 1\\
  \color{red}7&\color{red}3\\
  4\\
\end{array}\hspace{.25em} + \hspace{.25em} 
 \begin{array}{|c|c|}
  \color{red} 6 & 5\\
  2& \color{red}1 \\
  3&\color{red}4\\
  \color{red} 7 \\
\end{array}\hspace{.25em} + \hspace{.25em} 
\begin{array}{|c|c|}
  1 & 5\\
  \color{red} 6& \color{red}2 \\
  \color{red} 7&\color{red}3\\
  4 \\
\end{array}\hspace{.25em} + \hspace{.25em} 
\begin{array}{|c|c|}
  1 & 5\\
  \color{red} 6& \color{red}2 \\
  3&\color{red}4\\
  \color{red} 7 \\
\end{array}\hspace{.25em} + \hspace{.25em} 
\begin{array}{|c|c|}
  1 & 5\\
   2& \color{red}3\\
  \color{red}6&\color{red}4\\
  \color{red} 7 \\
\end{array}\hspace{.25em} 
 \hspace{.25em}   \right) .
 \end{align*}

\normalsize

In determining the coefficient of $[t]$ in Definition \ref{defeta}, we used the fact that ${m\choose \ell}$ is the number of ways to pick $\ell$ entries from the second column, making $\eta_\ell ([t])$ a sum of dual Garnir relations. 

Since $\eta_\ell$ is defined via its action on positions of $[t]$, it is a right action on $M^\lambda$. Meanwhile, $\sg_{n+m}$ acts on the letters of $t$, so its action on $\cal T_\lambda$ and its induced action on $M^\lambda$ are left actions. Therefore, the actions of $\eta_\ell$ and $\sg_{n+m}$ commute and  $\eta_\ell$ is $\sg_{n+m}$-equivariant. Furthermore, it follows from equation (\ref{Fultoneq}) that $\im(\eta_\ell) \subseteq \ker(\alpha)$, as $\eta_\ell([t])$ is a sum of dual Garnir relations. Using techniques employed in \cite{FHSW, FHW3, BF}, we will  show that for the vast majority of partitions $\mu$, the relations generated by $\eta_\ell$ for any single value of $\ell$ are all that is needed to generate ${G}^{\mu}$. To do so, we study the eigenvalues of $\eta_\ell$; the rest of this section is devoted to this study.

Note that because 
\[ {M}^{\mu} \cong \bigoplus_{i=0}^{m} S^{2^{i}1^{n+m - 2i}}\]
is multiplicity-free, by Schur's Lemma $\eta_\ell $ acts as a scalar on each irreducible submodule of ${M}^{\mu}$. Thus, finding the kernel of $\eta_\ell$ is equivalent to finding the irreducible submodules of ${M}^{\mu}$ on which $\eta_\ell $ acts like the 0 scalar.

We proceed by computing the action of $\eta_\ell$ on each irreducible submodule of ${M}^{\mu}$. For each $T \in \binom{[n+m]}{n}$, let $v_{T}\in {M}^{\mu}$ be the column tabloid with first column $T$ (both columns assumed to be in increasing order). For any $v\in {M}^{\mu}$, let $\langle v, v_T \rangle$ be the coefficient of $v_T$ in the expansion of $v$ in the basis of all $v_T$. 

\begin{lem}\label{lemma}
For every $S, T \in \binom{[n+m]}{n}$,  
\[ \langle \eta_\ell(v_{S}), v_{T} \rangle = \begin{cases}
{m \choose \ell} & \textrm{if } S = T, \\ 
0 & \textrm{if } S\neq T \; \; and \; \;  |S \cap T| \neq n-\ell  ,\\
(-1)^{\sum_{k=1}^\ell (c_k+d_k)+\ell +1} & \textrm{if } |S \cap T| = n -\ell  \textrm{ with}\\
&   S \backslash T=\{ c_1, \cdots , c_\ell \}, \; T \backslash S= \{ d_1, \ldots , d_\ell 
\}.

\end{cases} \]
\end{lem} 
\begin{proof}
The first two cases easily follow from the definition of $\eta _\ell$. We consider the third case.

Let $S=\{a_1, a_2, \ldots, a_n\}$ and $[n+m]\setminus S=\{b_1, b_2, \ldots, b_m\}$. Then 
 the columns of $v_S$, arranged from left to right instead of top to bottom, are
\[a_1, a_2, \ldots , a_n  \]
and
\[b_1, b_2,  \ldots , b_m,\]
where the entries of each column are in increasing order. Let $i_1, i_2, \ldots , i_\ell \in [n]$ and $j_1, j_2, \ldots , j_\ell \in [m]$. We will exchange the entries $a_{i_1}, a_{i_2}, \ldots , a_{i_\ell}$ with the entries $b_{j_1}, b_{j_2}, \ldots , b_{j_\ell}$ to obtain a term that appears in $\eta _\ell (v_S)$. The two columns of such a term look like
\[a_1, \ldots, a_{i_1-1}, {\color{red} b_{j_1}}, a_{i_1+1}, \ldots, a_{i_2-1}, {\color{red} b_{j_2}}, a_{i_2+1}, \ldots , a_{i_\ell -1}, {\color{red} b_{j_\ell}}, a_{i_\ell +1}, \ldots , a_n \]
and
\[b_1, \ldots, b_{j_1-1}, {\color{red} a_{i_1}}, b_{j_1+1}, \ldots, b_{j_2-1}, {\color{red}a_{i_2}}, b_{j_2+1}, \ldots , b_{j_\ell -1}, {\color{red}a_{i_\ell}}, b_{j_\ell +1}, \ldots, b_m.\]
In order to get an element in the basis of $M ^\mu$, we need to reorder the entries of each column so they are increasing. Every exchange within a column results in a sign. What sign will we end up with? We consider each column individually first. Each column will require two stages. Then we will put the resulting signs together. 

{\it Column 1, Stage 1}: We move all the $b$'s to the end of the column. Moving $b_{j_\ell}$ to the end gives $n-i_\ell$ transpositions. Moving $b_{j_{\ell -1}}$ to just before $b_{j_\ell}$ gives $n-i_{\ell -1}-1$ transpositions, and so on. The total number of transpositions for this stage is:
\begin{equation} \label{col1st1} \sum _{k=0}^{\ell -1} (n-i_{\ell -k}-k) = n\ell -\frac{\ell (\ell -1)}{2} - \sum _{k=1}^\ell i_k \; .\end{equation}
After these transpositions, the first column is
$$a_1, \ldots a_{i_1-1}, a_{i_1+1}, \ldots , a_{i_2-1}, a_{i_2+1}, \ldots a_{i_\ell -1}, a_{i_\ell +1}, \ldots, a_n, {\color{red} b_{j_1}}, {\color{red}b_{j_2}}, \ldots , {\color{red}b_{j_\ell}}\; .$$

{\it Column 1, Stage 2}: We move the $b$'s to their correct positions. We begin with $b_{j_1}$. We have that $b_{j_1}$ is larger than $j_1-1$ $b$'s and therefore it is larger than $(b_{j_1}-1)-(j_1-1)=b_{j_1}-j_1$ $a$'s. Some of those $a$'s, say $\gamma_1$ of them, are not in the first column (because they were exchanged into the second column). So in the first column, there are $b_{j_1}-j_1-\gamma _1$ $a$'s smaller than $b_{j_1}$.

We need to move $b_{j_1}$ from after all the $n-\ell$ $a$'s to after $b_{j_1}-j_1-\gamma _1$ $a$'s. That requires $(n-\ell)-(b_{j_1}-j_1-\gamma _1)$ transpositions. 

Define $\gamma_k$ to be the number of $a$'s in $\{a_{i_1}, \ldots , a_{i_\ell}\}$ that are smaller than $b_{j_k}$. So in the first column, there are $b_{j_k}-j_k-\gamma_k$ $a$'s smaller than $b_{j_k}$. 

We have dealt with $b_{j_1}$ already. Now we move $b_{j_2}$ from after all the $n-\ell$ $a$'s to after $b_{j_2}-j_2-\gamma_2$ $a$'s. Since $b_{j_1}<b_{j_2}$, we do not need to transpose through $b_{j_1}$. So we have $(n-\ell)-(b_{j_2}-j_2-\gamma _2)$ transpositions. Continuing this way gives 
\begin{equation} \label{col1st2} \sum _{k=1}^\ell \big ( (n-\ell)-(b_{j_k}-j_k-\gamma_k)\big )=\ell (n-\ell)-\sum_{k=1}^\ell \big ( b_{j_k}-j_k-\gamma _k \big ) \; \end{equation}
transpositions.
This completes the reordering of the first column. 

{\it Column 2, Stage 1}: We move all the $a$'s to the end of the column. By the same reasoning as in stage 1 of the first column, the number of transpositions required is 
\begin{equation}\label{col2st1} \sum _{k=0}^{\ell -1} (m-j_{\ell -k}-k) = m\ell -\frac{\ell (\ell -1)}{2} - \sum _{k=1}^\ell j_k \; .\end{equation}
Now the second column is
$$b_1, \ldots, b_{j_1-1}, b_{j_1+1}, \ldots , b_{j_2-1}, b_{j_2+1}, \ldots , b_{j_\ell -1}, b_{j_\ell +1}, \ldots, b_m, {\color{red}a_{i_1}}, {\color{red}\ldots a_{i_\ell}} \; .$$

{\it Column 2, Stage 2}: We  move the $a$'s to their correct positions. We have that $a_{i_k}$ is larger than $i_k-1$ of the $a$'s and $(a_{i_k} -1)-(i_k -1)=(a_{i_k}-i_k)$ of the $b$'s. 

Define $\delta _k$ to be the number of $b$'s in $\{b_{j_1}, \ldots , b_{j_\ell} \}$ that are smaller than $a_{i_k}$. Then we need to move $a_{i_1}$ left to just after $(a_{i_1}-i_1-\delta _1)$ $b$'s, which requires $(m-\ell)-(a_{i_k}-i_k-\delta _k)$ transpositions. Continuing with $a_{i_2}$ and so on, we have
\begin{equation} \label{col2st2} \sum_{k=1}^\ell \big ( (m-\ell)-(a_{i_k}-i_k-\delta _k) \big )= \ell (m-\ell)-\sum_{k=1}^\ell \big (a_{i_k}-i_k-\delta_k\big )\end{equation}
transpositions. 

{\it Total from both columns}: 
The total number of transpositions for both columns combined is given by the sum of equations (\ref{col1st1}), (\ref{col1st2}), (\ref{col2st1}), and (\ref{col2st2}). We are only concerned with the parity of the total number of transpositions. Adding the four equations and omitting any obviously even terms gives
\begin{equation*} \sum_{k=1}^\ell \big (\gamma _k+\delta _k \big )-\sum_{k=1}^\ell \big ( a_{i_k}+b_{j_k} \big )\end{equation*}
transpositions.

We now show that $\sum_{k=1}^\ell \big (\gamma _k+\delta _k \big )$ has the same parity as $\ell $. 

For any pair $(c,d)\in [\ell]\times [\ell]$, we have either $a_{i_c}<b_{j_d}$ or $a_{i_c}> b_{j_d}$. Therefore, we can think of the pair $(c,d)$ as contributing 1 to $\gamma_d$ in the first case, or as contributing 1 to $\delta_c$ in the second case. So each pair $(c,d)$ contributes exactly 1 to the sum $\sum_{k=1}^\ell \big (\gamma _k+\delta _k \big )$. Hence, the sum equals the number of pairs $(c,d)$, which is $\ell ^2$, which has the same parity as $\ell$.

Since $S\setminus T=\{a_{i_1}, \ldots , a_{i_\ell}\}$ and $T\setminus S=\{b_{j_1}, \ldots , b_{j_\ell}\}$, and there is a sign in the definition of $\eta _\ell$, the lemma follows. 
\end{proof}
We are now ready to compute the scalar action of $\eta_\ell$ on each irreducible submodule of $M^\mu$.
\begin{theorem} \label{etaell}
On the irreducible submodule of ${M}^{\mu}$ isomorphic to $S^{2^{i}1^{(n+m) - 2i}}$, the operator $\eta_\ell$ acts like multiplication by the scalar $\omega_{\ell, i}$, where
\small
\[\omega_{\ell, i}:= {m\choose \ell}-\sum_{\ell_1=0}^\ell {m-i\choose \ell_1}{n-i\choose \ell_1}{i\choose \ell -\ell_1}(-1)^{\ell_1}.\]
\normalsize
\end{theorem}
\begin{proof} Let $T=[n]$ so that
\small
$$v_T \, = \, 
\ytableausetup
{mathmode, boxsize=2.3em}
\begin{ytableau}
\scriptstyle1 &\scriptstyle n+1  \\
\scriptstyle 2 & \scriptstyle n+2  \\
\vdots & \vdots \\
\scriptstyle m &\scriptstyle n+m\\
\scriptstyle m + 1  \\
\vdots \\
\scriptstyle n \\
\end{ytableau}\hskip .5cm, 
$$
\normalsize
and let $t$ be the standard Young tableau of shape $2^i1^{n+m-2i}$ given by
\small
$$t\, = \, 
\ytableausetup
{mathmode, boxsize=2.3em}
\begin{ytableau}
\scriptstyle1 &\scriptstyle n+1  \\
\scriptstyle 2 & \scriptstyle n+2  \\
\vdots & \vdots \\
\scriptstyle i &\scriptstyle n+i\\
\scriptstyle i + 1  \\
\vdots \\
\scriptstyle n \\
\scriptstyle n+i+1\\
\vdots \\
\scriptstyle n+m 
\end{ytableau}\hskip .5cm.
$$
\normalsize
Recall that the Specht module $S^{2^{i}1^{(n+m) - 2i}}$is spanned by $\{\tau e_t :  \tau\in \sg_{n+m} \}$, 
where $e_t$ is the symmetrizer of Equation (\ref{et}). In order to study the action of $\eta_\ell$ on this Specht module, we begin by simplifying the action of $e_t$ on $v_T$ by factorizing $e_t$ as follows.  Let  $r_{t} = \sum_{\alpha \in R_{t}} \alpha$. Let $d_{t}$ be the signed sum of column permutations stabilizing $\{ 1, 2, \dots , n \}, \{ n+1, \dots, n+i \}$, and $\{ n+i+1, \dots n+ m \}$, i.e. the signed sum of permutations in the subgroup 
$$S_{\{ 1, \ldots , n\}}\times S_{\{n+1, \ldots , n+i\}} \times S_{\{n+i+1, \ldots , n+m\}} \subseteq C_t.$$ 
Now let $f_{t}$ be the signed sum of left coset representatives of the above subgroup of $C_t$, that is, permutations $\sigma$ in $C_{t}$ that satisfy\footnote{\label{leftrightfootnote}In \cite{FHSW} and \cite{BF}, the representatives used were mistakenly right coset representatives. The results in \cite{FHSW} are unaffected by the error. The main results in \cite{BF} are also unaffected, but see footnote on page \pageref{BFerror}. 
}
$$\sigma(1)<\cdots <\sigma(n), \hskip .3cm  \sigma(n+i+1)<\cdots <\sigma(n+m), \hskip .3cm \sigma(n+1)<\cdots <\sigma(n+i).$$
Then $e_{t}v_{T} = r_{t}f_{t}d_{t}v_{T}$. 
The antisymmetry of column tabloids ensures that $d_{t}v_{T}$ is a scalar multiple of $v_{T},$ because it simply permutes within columns. Therefore we can conclude that $r_{t}f_{t}v_{T}$ is a scalar multiple of $e_{t}v_{T}$, and in particular that $e_{t}v_{T}$ is nonzero, as the coefficient of $v_{T}$ in $r_{t}f_{t}v_{T}$ is 1.

Consider $\eta_\ell(r_{t}f_{t}v_{T}).$ In the subspace restricted to $S^{2^{i}1^{n+m-2i}}$, the fact that $\eta_\ell$ acts on $e_{t}v_{T}$ as a scalar implies the same is true of $r_{t}f_{t}v_{T}.$ In fact, because the coefficient of $v_{T}$ in $r_{t}f_{t}v_{T}$ is 1, we can determine precisely what this scalar is by computing $\langle \eta_\ell(r_{t}f_{t}v_{T}), v_{T} \rangle$.

We have 
\[ r_{t}f_{t}v_{T} = \sum_{S \in \binom{[n+m]}{n}} \langle r_{t}f_{t}v_{T}, v_{S} \rangle v_{S} .\]
Applying the linear operator $\eta_\ell$ thus gives 
\[ \eta_\ell(r_{t}f_{t}v_{T}) = \sum_{S \in \binom{[n+m]}{n}} \langle r_{t}f_{t}v_{T}, v_{S} \rangle \eta_\ell(v_{S}) .\]

Note that when $T = S$, by Lemma \ref{lemma} we have $\langle \eta_\ell(v_{T}), v_{T} \rangle = {m\choose \ell}$. With this, we can compute the coefficient of $v_{T}$ in general by 

\begin{align} \label{equation}
\omega_{\ell, i}= \langle \eta _\ell(r_{t} f_{t} v_{T}), v_{T} \rangle &=  \sum_{S \in \binom{[n+m]}{n}} \langle r_{t}f_{t}v_{T}, v_{S} \rangle \langle \eta_\ell(v_{S}), v_{T} \rangle  \\ \nonumber &= {m\choose \ell} + \sum_{S \in \binom{[n+m]}{n} \backslash \{ T \} }  \langle r_{t}f_{t}v_{T}, v_{S} \rangle \langle \eta_\ell (v_{S}), v_{T} \rangle .
\end{align}

The contributions to the sum arise only when $S$ and $T$ differ by $\ell$ elements. In the sum $r_{t}f_{t}v_{T}$, there are four different ways to obtain a $v_{S}$ that fulfills this criterion:
\begin{enumerate}
\item $v_S$ is obtained by doing  row permutations only (possible iff $\ell \leq i$).
\item $v_S$ is obtained by doing  column permutations only  (possible iff $\ell \leq m-i$). 
\item $v_S$ is obtained by a subset of the column permutations of case (2) above, followed by row permutations.
\item $v_S$ is obtained via a combination of  column permutations and  row permutations: $\ell_1$ entries are exchanged via column permutations and $\ell_2$ entries via row permutations, where $\ell=\ell_1+\ell_2$ (possible iff $\ell_1\leq m-i$ and $\ell _2\leq i$).
\end{enumerate}
In case (4), if we let $\ell _1=0$ we get case (1), and if we let $\ell_2=0$ we get cases (2) and (3). 
We find it helpful to compute cases (1), (2), and (3) before case (4). 
\vskip .2cm
\noindent {\bf {Case (1)}:  row exchanges only. }

Pick $u_1, u_2, \ldots , u_\ell \in [i]$ and let 
\[ \alpha = (u_1, n+u_1)(u_2, n+u_2) \cdots (u_\ell, n+u_\ell) \; .\]
Then the columns of $\alpha v_T$, written left to right instead of top to bottom, are
\[1, 2, \ldots , u_1-1, {\color{red} n+u_1}, u_1+1, \ldots, u_\ell -1, {\color{red} n+u_\ell}, u_\ell+1, \ldots , n\]
and
\[n+1, \ldots, n+u_1-1, {\color{red} u_1}, n+u_1+1, \ldots , n+u_\ell+1, {\color{red} u_\ell}, n+u_\ell+1, \ldots , n+m.\]
For $\alpha v_T$  to have ordered entries, we need to move the $n+u_1, n+u_2, \ldots , n+u_\ell$ to the end of the first column, and to
move the $u_1, u_2, \ldots , u_\ell$ in the second column to the beginning of the column. 

For the first column, we first move $n+u_\ell$ to the end, then $n+u_{\ell-1}$, and so on. The number of transpositions for the first column is:
\[ \sum_{k=0}^{\ell-1} (n-u_{\ell-k}-k)=n\ell-\sum_{k=1}^\ell u_k -\frac{\ell(\ell-1)}{2}.\]
The number of transpositions for the second column is:
\[ \sum_{k=1}^\ell (u_k-k) =\sum_{k=1}^\ell u_k -\frac{\ell(\ell+1)}{2}.\]
Together, the number of transpositions for ordering the columns of $\alpha v_T$ is
\[ n\ell -\ell ^2\; .\]

To compute $\langle \eta_\ell (v_S), v_T\rangle$, note that
$T\setminus S=\{ u_1, \ldots , u_\ell\}$ and $S\setminus T= \{n+u_1, \ldots , n+u_\ell\}$ in Lemma \ref{lemma}, giving 
\[  \langle \eta_\ell (v_S), v_T\rangle=(-1)^{\sum_{k=1}^\ell (u_k+n+u_k)+\ell +1}=(-1)^{n\ell +\ell +1}\; .\] 
Hence,
\[ \langle r_tf_tv_T, v_S\rangle \langle \eta_\ell (v_S), v_T\rangle=(-1)^{(n\ell +\ell +1)+(n\ell -\ell ^2) } =(-1) ,\]
which is independent of the choice of $u_k$'s. There are ${i\choose \ell}$ ways to pick the $u_k$'s. So case (1) gives an overall contribution of 
\begin{equation} -{i\choose \ell} .\end{equation}
 This expression appropriately gives 0 when $\ell >i$. 

\vskip .2cm
\noindent {\bf  Case (2): column exchanges only.}

The permutations $\sigma$ in $f_t$ satisfy
$$\sigma(1)<\cdots <\sigma(n), \hskip .3cm  \sigma(n+i+1)<\cdots <\sigma(n+m), \hskip .3cm \sigma(n+1)<\cdots <\sigma(n+i).$$ 
Since $\sigma \in C_t$, it follows that $\sigma$ fixes $\{n+1, \ldots , n+i \}$. Since we need $|S \cap T| = n-\ell$, and since $T=[n]$ and $S=\{\sigma(1), \ldots , \sigma(n)\}$, we need to exchange $\ell$ of the elements in $[ n]$ with $\ell$ of the elements $n+i+1, \ldots , n+m$. Let $Q_\ell := \{q_1< \ldots < q_\ell \}\in {[n]\choose \ell}$, and $P_\ell := \{ p_1< \ldots < p_\ell \}\in {\{n+i+1, \ldots , n+m\}\choose \ell}$ be the sets of elements exchanged. In one-line notation, $\sigma$ restricted to the first column of $t$ is the concatenation of the sequences 
$$1, \ldots, q_1-1,  q_1+1, \ldots, q_2-1, q_2+1, \ldots , q_\ell -1 ,  q_\ell +1, \ldots , n $$
$$ { p_1},{p_2}, \ldots , {p_\ell} $$
$${q_1}, {q_2}, \ldots , {q_\ell}$$ 
$$n+i+1, \ldots ,p_1-1,  p_1+1, \ldots , p_2-1, p_2+1, \ldots p_\ell -1,  p_\ell +1, \ldots , n+m.$$
The number of inversions in $\sigma$ restricted to the first column of $t$ is then
$$\sum _{k=0}^{\ell-1}(n-q_{\ell-k}-k)+\sum_{k=1}^\ell (p_k-k-n+i)+\ell^2 =i\ell +\sum_{k=1}^\ell (p_k-q_k).$$

Now we consider $\sigma v_T$. Its first column is the same as the first $n$ entries of the first column of $\sigma t$, and does not need further reordering. 
Its second column is the concatenation of the sequences
$$n+1, \ldots , n+i $$
$${q_1}, {q_2}, \ldots , {q_\ell}$$ 
$$n+i+1, \ldots ,p_1-1,  p_1+1, \ldots , p_2-1, p_2+1, \ldots p_\ell -1,  p_\ell +1, \ldots , n+m.$$
To reorder this  column, we need to move the $q_k$'s all the way to the left. That takes $i$ transpositions for each $q_k$, giving the sign $(-1)^{\ell i}$. So
$$\langle r_tf_tv_T,v_S\rangle = (-1)^{\sum_{k=1}^\ell (p_k-q_k)}.$$

It remains to compute $\langle \eta_\ell (v_S),v_T\rangle$. Since $S\setminus T=\{p_1, \ldots , p_\ell\}$ and $T\setminus S=\{q_1, \ldots , q_\ell\}$, in Lemma \ref{lemma} we have
$$\langle \eta_\ell (v_S),v_T\rangle = (-1)^{\sum_{k=1}^\ell (p_k+q_k)+\ell +1}.$$
The total sign of this contribution is then
$$(-1)^{\sum_{k=1}^\ell (p_k-q_k)}(-1)^{\sum_{k=1}^\ell (p_k+q_k)+\ell +1}=(-1)^{\ell +1}.$$
There are ${n\choose \ell}$ choices for the $q_k$'s and ${m-i\choose \ell}$ choices for the $p_k$'s, giving the total contribution from this case of
\begin{equation}  {n\choose \ell}{m-i\choose \ell}(-1)^{\ell +1}.\end{equation}
\vskip .5cm
\noindent {\bf  Case (3): column permutations from case (2) followed by row permutations.}

If in case (2), $q_j\leq i$ and $q_{j+1}>i$ for some $1\leq j\leq \ell$, then after $\sigma$ from case (2), we can apply row swaps $(q_\gamma, n+q_\gamma)$ for any $\gamma \leq j$ and still have $|S\cap T|=n-\ell$, with $S\setminus T=\{p_1, \ldots , p_\ell\}$ and $T\setminus S=\{q_1, \ldots , q_\ell\}$. The row group is not signed, so all that remains is  to reorder the columns of $\alpha \sigma v_T$ to be increasing and compute the corresponding sign.

Suppose we do the row swaps for  $q_{\gamma_1}, \ldots , q_{\gamma_r}$ for some $r\leq j$ and $1\leq \gamma_1 <\gamma_2<\cdots <\gamma_r\leq j$, so that
$$\alpha = (q_{\gamma_1}, n+q_{\gamma_1})\cdots (q_{\gamma_r}, n+q_{\gamma_r}).$$

The first column of $\alpha \sigma v_T$ is as in Case (2) and needs no further reordering. The sign for the reordering of second column of $\alpha \sigma v_T$ can be obtained from the sign of reordering of the second column of $\sigma v_T$ as follows. Suppose $\rho$ achieves the reordering of the second column of $\sigma v_T$. Then $\rho \alpha ^{-1}$ achieves the reordering of the second column of $\alpha \sigma v_T$. 
The contribution from case (2) is therefore modified only by $\sgn{\alpha^{-1}}=(-1)^r$, so we have in this case:
$$(-1)^{\ell +1+r}.$$
How many ways are there to have $r$ row swaps? In picking the $q_k$'s, let us pick $j$ of them to be in $[i]$ and $\ell -j$ of them to be in $[n]\setminus [i]$. Then we can pick $r$ of the first $j$ $q_k$'s for the row swaps. There are
\[{i\choose j}{n-i\choose \ell -j}{j\choose r}\]
ways to do this. There are still ${m-i\choose \ell}$ ways to pick the $p_k$'s. So we have the contribution
\begin{equation} {m-i\choose \ell}(-1)^{\ell +1} \sum _{r=1}^\ell \sum _{j=0}^\ell {i\choose j}{n-i\choose \ell -j}{j\choose r} (-1)^{r} .\end{equation}
Note that if we set $r=0$ in the above equation, we get the contribution of case (2).

\vskip .2cm
\noindent {\bf Case (4): column permutations to exchange $\ell _1$ elements and row permutations to exchange $\ell_2=\ell -\ell_1$ elements.}

We begin similarly to case (2), replacing $\ell$ with $\ell _1$, where $1\leq \ell_1<\ell$. We use a column permutation $\sigma$ to exchange the elements of $Q_{\ell_1}=\{q_1<\cdots <q_{\ell_1}\}\in {[n]\choose \ell_1}$ with the elements of $P_{\ell_1}=\{p_1<\cdots <p_{\ell_1}\}\in {\{n+i+1, \ldots , n+m\}\choose \ell_1}$. 
Then we pick $\ell_2=\ell-\ell_1$ elements of $[ i]$, say $1\leq u_1< \cdots < u_{\ell_2}\leq i$ such that $\{ u_1, \ldots , u_{\ell_2}\}\cap \{q_1, \ldots , q_{\ell_1}\}=\emptyset$ and carry out a row exchange on them,
\[\alpha = (u_1, n+u_1)(u_2, n+u_2) \cdots (u_{\ell_2}, n+u_{\ell_2}).\]
The permutation $\alpha \sigma$   exchanges a total of $\ell$ entries between the first and second columns of $v_T$, as desired. 

The sign of the column permutation is the same as in case (2), except that $\ell$ is replaced by $\ell _1$:
$$\sgn(\sigma)=(-1)^{-\ell_1 i +\sum_{k=1}^{\ell_1} (p_k-q_k)}.$$ 
We now consider the sign of the transposition that orders the two columns of $\alpha \sigma v_T$ in increasing order. We begin by applying only $\sigma$, so the first row of $\sigma v_T$ is
\[1, \ldots, q_1-1,  q_1+1, \ldots , q_{\ell_1} -1 ,  q_{\ell_1} +1, \ldots n, {\color{red} p_1},{\color{red}p_2}, \ldots , {\color{red}p_\ell},\]
and the second row of $\sigma v_T$ is
\[n+1, \ldots , n+i, {\color{red}q_1}, {\color{red}q_2}, \ldots , {\color{red}q_{\ell_1}}, n+i+1, \ldots ,p_1-1,  p_1+1, \ldots p_{\ell_1} -1,  p_{\ell_1} +1, \ldots n+m.\]

Once we apply $\alpha$, the entries $n+u_k$ are in the first column mixed in between 1 and $i$ where the $u_k$ used to be, and the entries $u_k$  are in the second column mixed in between $n+1$ and $n+i$ where the $n+u_k$ used to be. To count the number of transpositions required for reordering, we define $\gamma_k$ to be the number of elements in $Q_{\ell_1}$ larger than $u_k$. 

In the first column, we need to move the $n+u_k$'s to just after $n$. This requires
\[\sum_{k=1}^{\ell_2} (n-\gamma_k-u_k-(k-1))\]
transpositions. In the second column, we need to move the $u_k$'s and the $q_k$'s to the left and keep them in order. This requires
\[ i\ell_1+\sum_{k=1}^{\ell_1}(u_k-k+\gamma_k)\] 
transpositions. So turning $\alpha \sigma v_T$ into $v_S$ gives the sign
\[(-1)^{n\ell_2+i\ell_1+\ell_2}.\]
Finally, since $S\setminus T= \{p_1, \ldots , p_{\ell_1}, n+u_1, \ldots , n+u_{\ell_2}\}$ and $T\setminus S = \{p_1, \ldots, p_{\ell_1}, u_1, \ldots , u_{\ell_2}\}$, Lemma \ref{lemma} gives
\[ \langle \eta_\ell (v_S), v_T\rangle = (-1)^{ \sum_{k=1}^{\ell_1}(p_k+q_k)+\sum_{k=1}^{\ell_2}(n+u_k+u_k)+\ell+1}.\]
The overall sign simplifies to 
\[(-1)^{\ell_1+1}.\]
How many of these cases are there? If we choose $j$ of the $\ell_1$ $q$'s to come from $[i]$, there will be $i-j$ elements in $[i]$ to choose the $u_k$'s from. Summing these and remembering that there are still ${m-i\choose \ell_1}$ ways to choose the $p_k$'s gives
\[ {m-i\choose \ell_1}(-1)^{\ell_1+1}\sum _{j=0}^{\ell_1} {i\choose j}{n-i\choose \ell_1-j}{i-j\choose \ell_2}.\]
As in case (3), we can also have row swaps using the subset of the $q_k$'s that are in $[i]$, and if there are $r$ such row swaps, this merely introduces a factor of $(-1)^r$, giving
\begin{equation} {m-i\choose \ell_1}(-1)^{\ell_1+1}\sum_{r=0}^j\sum _{j=0}^{\ell_1} {i\choose j}{n-i\choose \ell_1-j}{i-j\choose \ell_2}{j\choose r}(-1)^{r}.\end{equation}
We start the sum from $r=0$ to include the case we just computed where no additional row swaps are done. 

Note that setting $\ell_1=0$ in the above formula gives us the contribution from case (1) and setting $\ell_1=\ell$ in the above formula gives us the contributions from cases (2) and (3) combined. Therefore, the total eigenvalue from all four cases is given by summing the above formula from $\ell_1=0$ to $\ell_1=\ell$.

Now, reordering the sum to carry out the sum over $r$ first reveals that we have a term
\begin{equation} \label{cancellation}\sum_{r=0}^j{j\choose r}(-1)^r = 
\begin{cases} 1&j=0\\ 0& j>0 \end{cases}.\end{equation}
So the only contribution to the sum above is the $j=r=0$ term, and the theorem is proved. 

Note that the cancellation that appears in equation (\ref{cancellation}) means in essence that we can limit the $Q_{\ell_1}$ to only ${[n]\setminus [i] \choose \ell_1}$ and not consider  the  contributions from the type of row swaps discussed in case (3) (where the $r$ originates), since those cancel the contributions from the $Q_{\ell_1}$ in ${[n]\choose \ell_1}\setminus {[n]\setminus [i] \choose \ell_1}$. 

\end{proof}

\noindent {\bf Example}. Let $n=8$, $m=7$, $\ell=3$. We will demonstrate case (4) with  $\ell_1=2$, $\ell_2=1$, $i=4$. 

Pick $q_1=2$, $q_2=5$, $p_1=14$, $p_2=15$, $u_1=3$. In this case, $j=1$ (because $q_1\leq 4$ and $q_2>4$). We have: 
 
$$v_T\, = \, 
\ytableausetup
{mathmode, boxsize=1.7em}
\begin{ytableau}
\scriptstyle1 &\scriptstyle 9  \\
\scriptstyle 2 & \scriptstyle 10  \\
\scriptstyle 3 &\scriptstyle 11  \\
\scriptstyle 4 & \scriptstyle 12  \\
\scriptstyle 5 &\scriptstyle 13  \\
\scriptstyle 6 & \scriptstyle 14  \\
\scriptstyle 7 &\scriptstyle 15  \\
\scriptstyle 8    \\
\end{ytableau}\hskip .2cm ,
\hskip .2cm t\, = \, 
\ytableausetup
{mathmode, boxsize=1.7em}
\begin{ytableau}
\scriptstyle 1 &\scriptstyle 9  \\
\scriptstyle 2 & \scriptstyle 10  \\
\scriptstyle 3 &\scriptstyle 11  \\
\scriptstyle 4 & \scriptstyle 12  \\
\scriptstyle 5 \\
\scriptstyle 6  \\
\scriptstyle 7 \\
\scriptstyle 8\\
\scriptstyle 13 \\
\scriptstyle 14\\
\scriptstyle 15\\
\end{ytableau}\hskip .2cm,\hskip .2cm 
\sigma t\, = \, 
\ytableausetup
{mathmode, boxsize=1.7em}
\begin{ytableau}
\scriptstyle 1 &\scriptstyle 9  \\
\scriptstyle 3 & \scriptstyle 10  \\
\scriptstyle 4 &\scriptstyle 11  \\
\scriptstyle   6 & \scriptstyle 12  \\
\scriptstyle 7 \\
\scriptstyle  8  \\
\scriptstyle \color{red}14 \\
\scriptstyle \color{red}15\\
\scriptstyle  \color{red}\color{red}2 \\
\scriptstyle \color{red}5\\
\scriptstyle  13 \\
\end{ytableau}\hskip .2cm, \hskip .2cm 
\alpha \sigma v_T\, = \, 
\ytableausetup
{mathmode, boxsize=1.7em}
\begin{ytableau}
\scriptstyle1 &\scriptstyle 9  \\
\scriptstyle \color{blue} 11 & \scriptstyle 10  \\
\scriptstyle 4 &\scriptstyle \color{blue} 3 \\
\scriptstyle  6 & \scriptstyle 12  \\
\scriptstyle 7&\scriptstyle   \color{red}2  \\
\scriptstyle  8 & \scriptstyle  \color{red}5  \\
\scriptstyle  \color{red}14 &\scriptstyle  13  \\
\scriptstyle  \color{red}15    \\
\end{ytableau}\hskip .2cm , \hskip .2cm
v_S\, = \, 
\ytableausetup
{mathmode, boxsize=1.7em}
\begin{ytableau}
\scriptstyle1 &\scriptstyle \color{red} 2  \\
\scriptstyle 4 & \scriptstyle \color{blue} 3  \\
\scriptstyle 6 &\scriptstyle \color{red}5 \\
\scriptstyle  7 & \scriptstyle 9 \\
\scriptstyle 8&\scriptstyle   10 \\
\scriptstyle  \color{blue} 11& \scriptstyle 12  \\
\scriptstyle  \color{red}14 &\scriptstyle  13  \\
\scriptstyle  \color{red}15    \\
\end{ytableau}\hskip .2cm .
$$
Following the counting in the proof for $\sgn (\sigma)$, we have
$$\sgn(\sigma)=(-1)^{-2\cdot 4+(14-2)+(15-5)}=+1.$$
Alternatively, we see that in cycle notation, $\sigma = (2,3,4,6,8,15,13)(5,7,14)$, which is even. 

Now following the counting in the proof for ordering the columns of $\alpha \sigma v_T$, we have
$$(-1)^{8\cdot 1+4\cdot 2+1}=(-1).$$
Alternatively, to order the first column we can use the permutation $(11,4,6,7,8)$, which is even, and to order the second column we can use the permutation $(9,2, 10, 3,5,12)$, which is odd, giving the overall sign $(-1)$.

Finally, since $T\setminus S=\{2,3,5\}$ and $S\setminus T=\{11,14,15\}$, the sign we get from Lemma \ref{lemma} for $\langle \eta_2(v_S),v_T\rangle$ is 
$$(-1)^{(2+3+5)+(11+14+15)+3+1}=+1.$$

So the contribution of this case to the eigenvalue is $$(+1)(-1)(+1))=-1=(-1)^{\ell_1+1},$$
where $\ell_1=2$.
\vskip .3cm
\noindent {\it The trace:} The operator  $\eta_2$ for the shape $(8,7)'$ is a ${15\choose 8}\times {15\choose 8}$ matrix. Its diagonal entries all equal ${7\choose 3}$, so its trace  is ${7\choose 3}{15\choose 8}$. A computation of the sum of the eigenvalues of $\eta_2$ with multiplicity, where the multiplicity of $\omega_{2,i}$ is the dimension of $S^{2^i1^{15-2i}}$ for $i=0, 1, \ldots , 7$, confirms this trace. 
The computation of the trace in these two different ways for any $n\geq m\geq \ell$ gives the following result. 

\begin{corollary}[A combinatorial identity] \label{id}For $n\geq m\geq \ell$, the following identity holds:
\[ \sum_{i=0}^m \sum_{\ell_1=0}^\ell {m-i\choose \ell_1}{n-i\choose \ell_1}{i\choose \ell -\ell_1} {n+m\choose i}\frac{n+m-2i+1}{n+m-i+1}(-1)^{\ell_1}=0.\]
\end{corollary}
\begin{proof} Given $n$, $m$, and $\ell$, the diagonal entries of $\eta_\ell$ are equal ${m\choose \ell}$, and its size is $\dim M^{(n,m)'}\times \dim M^{(n,m)'}$. Since $\dim  M^{(n,m)'}= {n+m\choose n}$, we have
$$\mbox{Tr } (\eta_\ell) = {m\choose \ell}{n+m\choose n}.$$
The trace also equals the sum of the eigenvalues of $\eta_\ell$ with multiplicity, i.e. 
\[ \mbox{Tr }(\eta_\ell)=\sum_{i=0}^m \omega_{\ell, i}\dim (S^{2^i1^{n+m-2i}}).\]
Hooke's Law formula gives
\[\dim S^{2^i1^{n+m-2i}}= {n+m\choose i}\frac{n+m-2i+1}{n+m-i+1}.\]
So we have
\small
\[ {m\choose \ell}{n+m\choose n}=\sum_{i=0}^m \left ({m\choose \ell}-\sum_{\ell_1=0}^\ell {m-i\choose \ell_1}{n-i\choose \ell_1}{i\choose \ell -\ell_1}(-1)^{\ell_1}\right) {n+m\choose i}\frac{n+m-2i+1}{n+m-i+1}.\]
\normalsize
Since
\[ {M}^{(n,m)'} \cong \bigoplus_{i=0}^{m} S^{2^{i}1^{n+m - 2i}},\]
we have the identity
\[{n+m\choose n}=\sum_{i=0}^m \dim (S^{2^{i}1^{n+m - 2i}}) = \sum_{i=0}^m {n+m\choose i}\frac{n+m-2i+1}{n+m-i+1},\]
which results in the corollary. 
\end{proof}

\section{A class of presentations of Specht modules}\label{class}

It is straightforward to check that  for all 2-column partitions $\mu = 2^m1^{n-m}$, we have $w_{\ell,m}=0$. It follows that
\[\ker \eta_\ell \supseteq S^{\mu}.\]
The Specht modules for which $w_{\ell,i}\neq 0$ for $0\leq i< m$ are those for which we have obtained a new presentation. That is, in those cases, $\ker(\eta_\ell) \cong S^\mu$ and $\im(\eta_\ell) = \ker(\alpha)$ for $\alpha : {M}^\mu \rightarrow S^\mu$. This is the content of our central theorem. 

\begin{theorem} \label{imetaellkeralpha}
Let $\mu = 2^m1^{n-m}$ and let $H^{\mu, \ell}$ be the subspace of $M^\mu$ generated by $\eta_\ell([t])$ for $[t]\in M^\mu$. Then
\[ M^\mu/H^{\mu,\ell}  \cong  S^\mu\] as $\sg_{n+m}$-modules iff
\small
\[ {m\choose \ell}-\sum_{\ell_1=0}^\ell {m-i\choose \ell_1}{n-i\choose \ell_1}{i\choose \ell -\ell_1}(-1)^{\ell_1}\neq 0 \]
\normalsize
for $i=0, 1, \ldots, m-1$.

\end{theorem}
\begin{proof} The theorem is a direct consequence of Theorem \ref{etaell}. \end{proof}

We now state the conditions for presentations of partitions with two or more columns. Let $\lambda$ be a partition, let $\ell_c\in [\lambda'_{c+1}]$, and let $h_{c,\ell_c}([t])$ be the image of $\eta_{\ell_c}$ on the $c$ and $(c+1)^{st}$ columns of $[t]$ that leaves the other columns of $[t]$ fixed.  

\begin{theorem}\label{spechtgarnirell}
Let $\lambda$ be a partition of $n$ and let $\lambda'$ be its conjugate partition. Let $\hat \ell=(\ell_1, \ell _2, \ldots, \ell_{\lambda_1-1})$ and let $H^{\lambda, \hat \ell}$ be the space generated by $h_{c,\ell_c}([t])$ for every $[t]\in M^\lambda$, $1\leq c \leq \lambda _1-1$, and one choice of $\ell_c\in [\lambda'_{c+1}] $ for each $c$. Then 
\[       M^\lambda/H^{\lambda, \hat\ell} \cong_{\sg_n} S^\lambda \]
iff $\lambda$ satisfies the conditions 
\[ {\lambda '_{c+1}\choose \ell_c}-\sum_{k=0}^{\ell_c}{\lambda '_{c+1}-i\choose k}{\lambda '_c-i\choose k}{i\choose \ell _c -k}(-1)^{k}\neq 0 \]
\normalsize
for all $1\leq c\leq \lambda_c-1$ and $i=0,1,\ldots , \lambda'_{c+1}-1$. 

\end{theorem}

\begin{proof} This theorem follows from Theorem \ref{imetaellkeralpha} and the definition of $h_{c, \ell_c}([t])$. \end{proof} 

\noindent {\bf Remark} \label{MMSrelation} Motivated by the same question proposed in \cite{FHW3} that led to the present paper, \cite{MMS} addressed the case of intermediate values of $\ell$, approaching it via representations of the general linear group. In that different context, but using the same generalization of the symmetrized sum of \cite{BF} (Equation (\ref{BFsym})) as we do in Definition \ref{defeta}, they state the following sufficient (but not necessary, except in the 2-column case) condition for the corresponding quotient space to be a Specht module. We modify their notation to match ours.

\begin{theorem}\label{MMS}\cite[Theorem 6.2]{MMS}\footnote{Due to a typo, the statement of this theorem in the published version of \cite{MMS} erroneously limits the values of $j$ to $j=1, \ldots , \ell_c$; the values should be $j=1, \ldots ,\lambda_{c+1}'$, as stated here.} Let $\lambda$ be a partition of $n$ and $\ell_1, \ldots , \ell_{\lambda_1-1}$ positive integers satisfying $\ell_c\leq\lambda_{c+1}$ , $c= 1, \ldots , \lambda_1 -1$. Let $\hat \ell = (\ell_1, \ldots , \ell_{\lambda_1 -1})$. Then as $\sg_n$-modules, we have $M^\lambda / H^{\lambda, \hat \ell}\cong S^\lambda$ if
\[\sum_{t=1}^j (-1)^{t-1}{\lambda_{c+1}'-t\choose \lambda_{c+1}'-\ell_c}{j\choose t}{\lambda_c'-\lambda _{c+1}'+j+t\choose t} \neq 0\]
for all $c=1, \ldots , \lambda_1-1$ and $j=1, \ldots , \lambda_{c+1}'.$
\end{theorem}
For the two-column case of Theorem \ref{MMS}, i.e. for $\lambda_1=2$, \cite{MMS} state that this condition is also necessary (\cite[Corollary 6.1]{MMS}). Putting their result together with Theorem \ref{imetaellkeralpha} gives

\begin{corollary}[Equivalence of two sets of combinatorial conditions]\label{equivcomb} Let $n\geq m\geq \ell$. Then 
\[ {m\choose \ell}-\sum_{\ell_1=0}^\ell {m-i\choose \ell_1}{n-i\choose \ell_1}{i\choose \ell -\ell_1}(-1)^{\ell_1}\neq 0 \]
for $i=0, 1, \ldots, m-1$ 
iff 
\[ \sum_{t=1}^j (-1)^{t-1} {m-t\choose m-\ell}{j\choose t}{n-m+j+t\choose t}\neq 0 \]
for $j=1, \ldots , m$.
\end{corollary}
\vskip .2cm

As in \cite{BF}, Theorem \ref{spechtgarnirell} dramatically reduces the number of generators needed to obtain ${G}^{\lambda}$. The original construction leading to Equation (\ref{Fultoneq}) required enumerating over every $1 \leq k \leq \lambda _{c+1}'$ for every pair of columns $c$ and $c+1$ of every $t \in \Tlam$. Even Fulton's simplification using only $g_{c,1}$ relations requires enumerating over $t \in \Tlam$ for every pair of columns $c$ and $c+1$. By contrast, our construction uses a single relation for every pair of adjacent columns, and $[t]$ varies in $M^\lambda$, a significantly smaller space than $\Tlam$.

We now discuss the shapes $\mu$ and values $\ell$ for which the conditions of Theorem \ref{imetaellkeralpha} hold. 

When $\ell=1$, the condition $w_{1,i}\neq 0$ in Theorem \ref{imetaellkeralpha} simplifies to $(m-i)(n-i+1)\neq 0$, which is achieved in all cases, i.e. whenever $0\leq i <m\leq n$. This is the case of \cite{BF},\footnote{\label{BFerror}  In \cite{BF}, the error described in the footnote on page \pageref{leftrightfootnote} led to the incorrect formula $(m-i)(n+1)\neq 0$ which, like the correct formula obtained above, is also achieved in all cases; that error has been corrected here, and Theorem \ref{bf}, the main result of \cite{BF}, holds.} cited earlier as Theorem \ref{bf}.

When $\ell=m$, the condition $w_{m,i}\neq 0$ in Theorem \ref{imetaellkeralpha} simplifies to 
$$1-{n-i \choose m-i}(-1)^{m-i}\neq 0,$$
which holds for $0\leq i <m$ whenever $n\neq m$ or $n=1$. This is the case of \cite{FHW3}, cited earlier as Theorem \ref{fhw3}. In fact, when $\ell=m$, no symmetrization occurs and $\eta_m$ is equal to the Garnir relation $g_{1,m}$ used in \cite{FHW3}.

Generally speaking, the conditions in Theorem \ref{imetaellkeralpha} seem to hold for the vast majority of cases. Considering $n\geq m \geq \ell$ for $1\leq n \leq 50$, only 391 of the possible 22,100 combinations of $n$, $m$, and $\ell$ have values of $i<m$ for which $S^{2^{i}1^{n+m - 2i}}$ is in the kernel of $\eta_\ell$. In those cases,  there is only one such value of $i$ in 
all but 12 of the cases, in which the number of such values is 2. Data for $1\leq n\leq 28$ is recorded in Table \ref{tab:wli}.

\counterwithin{table}{section}
\begin{table}[h] \caption{Values of $(n,m,\ell,i)$ for which $w_{\ell, i}=0$ for $n\leq 28$, $0\leq i<m$, and $1\leq \ell<m$. }  \label{tab:wli}
\footnotesize
\begin{tabular}{|r|r|r|r|} 
\hline 
 $n$&$m$& $\ell$ & $i$ \\
\hline \hline
5 & 4 & 2 & 1 \\
6 & 4 & 3 & 2 \\
6 & 6 & 4 & 3 \\
6 & 6 & 5 & 3 \\
7 & 3 & 2 & 1 \\
7 & 5 & 3 & 1 \\
7 & 6 & 5 & 3 \\
7 & 6 & 5 & 4 \\
7 & 7 & 5 & 4 \\
8 & 5 & 2 & 2 \\
8 & 5 & 4 & 2 \\
8 & 8 & 2 & 3 \\
9 & 5 & 3 & 3 \\
9 & 7 & 2 & 3 \\
9 & 7 & 3 & 2 \\
9 & 7 & 5 & 5 \\
9 & 9 & 4 & 4 \\
9 & 9 & 4 & 5 \\
9 & 9 & 5 & 5 \\
10 & 9 & 2 & 4 \\
10 & 9 & 5 & 5 \\
10 & 9 & 5 & 6 \\
10 & 9 & 7 & 6 \\
10 & 10 & 5 & 6 \\
11 & 6 & 2 & 3 \\
11 & 7 & 4 & 5 \\
11 & 7 & 5 & 4 \\
11 & 8 & 5 & 6 \\
11 & 11 & 2 & 5 \\
11 & 11 & 9 & 5 \\
12 & 4 & 2 & 2 \\
12 & 6 & 3 & 2 \\
12 & 6 & 3 & 4 \\
12 & 10 & 2 & 5 \\
12 & 10 & 3 & 6 \\
12 & 10 & 7 & 8 \\
12 & 10 & 9 & 6 \\
12 & 11 & 9 & 9 \\
12 & 12 & 4 & 6 \\
12 & 12 & 8 & 6 \\
\hline
\end{tabular}
\hskip .1cm
\begin{tabular}{|r|r|r|r|} 
\hline 
 $n$&$m$& $\ell$ & $i$ \\
\hline \hline
13 & 6 & 4 & 3 \\
13 & 9 & 5 & 7 \\
13 & 10 & 4 & 7 \\
13 & 10 & 7 & 7 \\
13 & 12 & 7 & 9 \\
14 & 7 & 2 & 4 \\
14 & 11 & 2 & 6 \\
14 & 12 & 3 & 5 \\
14 & 12 & 6 & 8 \\
15 & 7 & 3 & 5 \\
15 & 10 & 5 & 8 \\
15 & 11 & 6 & 9 \\
15 & 13 & 9 & 11 \\
15 & 14 & 10 & 11 \\
15 & 14 & 11 & 8 \\
15 & 15 & 13 & 10 \\
15 & 15 & 14 & 10 \\
16 & 9 & 3 & 4 \\
16 & 10 & 2 & 6 \\
16 & 12 & 2 & 7 \\
16 & 14 & 3 & 9 \\
16 & 15 & 7 & 11 \\
16 & 15 & 14 & 10 \\
16 & 15 & 14 & 11 \\
16 & 16 & 4 & 5 \\
16 & 16 & 14 & 11 \\
17 & 5 & 2 & 3 \\
17 & 8 & 2 & 5 \\
17 & 11 & 5 & 9 \\
17 & 13 & 7 & 11 \\
17 & 15 & 3 & 7 \\
17 & 15 & 13 & 7 \\
17 & 16 & 2 & 9 \\
17 & 16 & 10 & 13 \\
17 & 16 & 11 & 13 \\
17 & 16 & 13 & 14 \\
18 & 8 & 3 & 6 \\
18 & 9 & 5 & 6 \\
18 & 10 & 4 & 6 \\
18 & 10 & 4 & 8 \\
\hline
\end{tabular}
\hskip .2cm
\begin{tabular}{|r|r|r|r|} 
\hline 
 $n$&$m$& $\ell$ & $i$ \\
\hline \hline
18 & 12 & 7 & 9 \\
18 & 13 & 2 & 8 \\
18 & 14 & 5 & 11 \\
18 & 15 & 9 & 13 \\
18 & 16 & 11 & 13 \\
18 & 16 & 11 & 14 \\
18 & 16 & 14 & 12 \\
18 & 17 & 11 & 10 \\
18 & 17 & 11 & 14 \\
19 & 9 & 6 & 3 \\
19 & 12 & 5 & 3 \\
19 & 12 & 5 & 10 \\
19 & 15 & 8 & 13 \\
20 & 9 & 2 & 6 \\
20 & 14 & 2 & 9 \\
20 & 16 & 2 & 10 \\
20 & 17 & 7 & 14 \\
20 & 17 & 11 & 14 \\
21 & 9 & 3 & 7 \\
21 & 13 & 5 & 11 \\
21 & 15 & 3 & 11 \\
21 & 17 & 9 & 15 \\
21 & 19 & 13 & 17 \\
22 & 6 & 2 & 4 \\
22 & 15 & 2 & 10 \\
22 & 15 & 4 & 12 \\
22 & 15 & 8 & 12 \\
22 & 16 & 7 & 14 \\
22 & 19 & 2 & 12 \\
22 & 21 & 17 & 19 \\
23 & 10 & 2 & 7 \\
23 & 12 & 9 & 7 \\
23 & 13 & 2 & 9 \\
23 & 14 & 5 & 12 \\
23 & 16 & 11 & 12 \\
23 & 19 & 10 & 17 \\
23 & 21 & 9 & 17 \\
23 & 23 & 2 & 14 \\
24 & 10 & 3 & 8 \\
24 & 15 & 4 & 10 \\
\hline
\end{tabular}
\hskip .2cm
\begin{tabular}{|r|r|r|r|} 
\hline 
 $n$&$m$& $\ell$ & $i$ \\
\hline \hline
24 & 16 & 2 & 11 \\
24 & 16 & 6 & 14 \\
24 & 19 & 9 & 17 \\
24 & 22 & 2 & 14 \\
24 & 22 & 3 & 12 \\
24 & 22 & 11 & 19 \\
24 & 22 & 15 & 20 \\
24 & 22 & 20 & 12 \\
24 & 24 & 4 & 15 \\
25 & 13 & 4 & 11 \\
25 & 15 & 5 & 11 \\
25 & 15 & 5 & 13 \\
25 & 21 & 11 & 19 \\
25 & 22 & 13 & 20 \\
25 & 23 & 9 & 19 \\
26 & 11 & 2 & 8 \\
26 & 17 & 2 & 12 \\
26 & 17 & 3 & 10 \\
26 & 22 & 8 & 19 \\
26 & 22 & 13 & 19 \\
26 & 24 & 3 & 17 \\ 
26 & 25 & 2 & 16 \\
27 & 7 & 2 & 5 \\
27 & 11 & 3 & 9 \\
27 & 16 & 5 & 14 \\
27 & 19 & 7 & 17 \\
27 & 21 & 9 & 19 \\
27 & 22 & 2 & 15 \\
27 & 23 & 12 & 21 \\
27 & 25 & 17 & 23 \\
27 & 26 & 21 & 24 \\
28 & 18 & 2 & 13 \\
28 & 26 & 3 & 15 \\
28 & 28 & 2 & 18 \\
28 & 28 & 26 & 18 \\
28 & 28 & 26 & 21 \\
28 & 28 & 27 & 21 \\
&&&\\
&&&\\
&&&\\
\hline
\end{tabular}
\normalsize
\end{table}

\vskip .2cm

\vskip .5cm
\noindent {\bf Acknowledgements}

 The author is grateful to Phil Hanlon and Michelle Wachs  for related collaboration and helpful discussions; to Sarah Brauner for helpful discussions and comments on an earlier draft;  to Thomas McElmurry for helpful discussions; and to Tony Yan for writing the computer program for calculating $w_{\ell , i}$.

\end{document}